\documentclass[11pt]{article}
\usepackage{amssymb, amsmath,  amsfonts,dsfont, mathrsfs,euscript,eufrak,colonequals}
\usepackage{textcomp}
\usepackage{latexsym}
\usepackage{indentfirst}
\usepackage[normalem]{ulem}

\usepackage{graphicx}
\usepackage[colorlinks=true, allcolors=blue]{hyperref}

\newtheorem{Theorem}{Theorem}

\newcommand{\N}{\mathbb N}
\newcommand{\R}{\mathbb R}
\newcommand{\Z}{\mathbb Z}

\newcommand{\CC}{\mathbb C}

\newcommand{\oo}{\bar{0}}

\newcommand{\lng}{\langle}
\newcommand{\rng}{\rangle}

\newcommand{\dd}{d}
\newcommand{\Ln}{\mathrm{Ln}\,}

\newcommand{\e}{\varepsilon}

\newcommand{\defeq}{\colonequals}

\newcommand{\Index}{\mathcal{I}} 

\begin{document}
	\author{I. A. Alexeev\footnote{St. Petersburg Department 
		of Steklov Mathematical Institute 
		of Russian Academy of Sciences, 27 Fontanka, Saint-Petersburg} \footnote{Institute for Information Transmission Problems of Russian Academy of Sciences, Bolshoy Karetny per. 19, build.1, Moscow, email: \texttt{vanyalexeev@list.ru}}, A. A. Khartov\footnote{Laboratory for Approximation Problems of Probability,  Smolensk State University, 4 Przhevalsky st., 214000 Smolensk, Russia, web-site: \texttt{approlab.org}} \footnote{Scientific and Educational Center of Mathematics, ITMO University,	Kronverksky Pr. 49,	197101, Saint-Petersburg, Russia, e-mail: \texttt{alexeykhartov@gmail.com}}}

\title{A criterion and a Cram\'er--Wold device for quasi-infinite divisibility  for discrete multivariate probability laws}
\date{}

\maketitle
\begin{abstract}
	Multivariate discrete probability laws  are considered. We show that such laws are quasi-infinitely divisible if and only if their characteristic functions are separated from zero. We generalize  the existing results for the univariate discrete laws and for the multivariate laws on $\Z^d$. The Cram\'er--Wold devices for infinite and quasi-infinite divisibility were proved. 
\end{abstract}


\section{Introduction. }

Let $F$ be a  distribution function of a multivariate probability law on $\R^d$, where $\R$ is the real line, $d$ is a positive integer. Recall that $F$ and the corresponding law are called \textit{infinitely divisible} if for every positive integer $n$ there exists a distribution function $F_n$ such that $F = F_n^{*n}$, where ``$*$"\, denotes the convolution, i.e. $F$ is the $n$-fold convolution power of $F_n$.  It is known that $F$ is infinitely divisible if and only if its characteristic function 
\begin{eqnarray*}
	f(t) \defeq \int_{\R^d} e^{i\lng t,x\rng} \dd F(x), \quad t \in \R^d,
\end{eqnarray*}
admits the following \textit{L\'evy representation} (see \cite[Theorem 8.1]{Sato_book})
\begin{eqnarray}\label{ID_Levy_representation}
	f(t) = \exp \biggl\{ i\lng t,\gamma\rng - \tfrac{1}{2}\lng t, Qt\rng + \int_{\R^d} \Bigr(e^{i\lng t,x\rng} - 1 - \tfrac{i\lng t,x\rng}{1+\|x\|^2} \Bigr) \nu(\dd x) \biggr\},\quad t\in \R^d,
\end{eqnarray}
where $\lng\, \cdot\,, \cdot\, \rng$ denotes the standard scalar product in $\R^d$, $\|x\|\defeq \sqrt{\lng x,x\rng}$ for any $x\in\R^d$,  $\gamma \in \R^d$ is a fixed vector, $Q$ is a symmetric nonnegative-definite $d\times d$ matrix, and $\nu$ is a measure on $\R^d$ that satisfies the following conditions
\begin{eqnarray*}
	\nu\bigr( \{\oo\}\bigr) = 0, \quad \int_{\R^d} \min\bigr\{\|x\|^2,1\bigr\} \nu(\dd x) < \infty. 
\end{eqnarray*}
Here and below, we denote by $\oo$  the zero vector of $\R^d$. The vector $(\gamma, Q, \nu)$ is called a\textit{ characteristic triplet} and it is uniquely determined	by $f$ and hence by $F$.

In the recent paper by Berger, Kutlu, and Lindner  \cite{BergKutLind}, the notion of quasi-infinitely divisible distributions on $\R^d$ was introduced. Following them, a distribution function $F$ and the corresponding law are called \textit{quasi-infinitely divisible}, if there exist infinitely divisible distribution functions $F_1$ and $F_2$ such that $ F_1=F * F_2$ (in the papers \cite{Khartov2}--\cite{Khartov4}, such property is proposed to be called \textit{rational infinite divisibility}). It was proved in \cite{BergKutLind} that $F$ is quasi-infinitely divisible if and only if the representation \eqref{ID_Levy_representation}  holds, where $\nu$ is a a signed finite measure on $\R^d \setminus (-r, r)^d$ for any $r>0$ that satisfies $	\nu\bigr( \{\oo\}\bigr) = 0$, and 
\begin{eqnarray*}
	\int_{\R^d} \min\bigr\{\|x\|^2,1\bigr\} |\nu|(\dd x) < \infty,
\end{eqnarray*}
where $|\nu|$ denotes the total variation of the measure $\nu$ (see \cite{BergKutLind} for more details). It is seen that quasi-infinitely divisible distributions are natural generalizations of infinitely divisible distributions. 

The examples of  univariate quasi-infinitely divisible laws can be found in the classical monographs \cite{GnedKolm}, \cite{LinOstr}, and  \cite{Lukacs}. The first detailed analysis of these laws on $\R$ was performed in \cite{Lindner}, and a lot of results for the univariate case are contained in the works \cite{Alexeev_Khartov}, \cite{Berger}, \cite{BergKut}, \cite{Khartov1}, \cite{Khartov2}, and \cite{Khartov3}. The multivariate case was considered in the recent papers \cite{BergKutLind}, \cite{Lindner_Berger}, and \cite{Passeggeri}. In these works the authors studied  questions concerning  supports, moments, continuity, and the weak convergence. The most complete results were obtained for probability laws on the set $\Z^d$, where $\Z$  is the set of integers. In particular, the following important fact was stated in \cite{Lindner_Berger}.
\begin{Theorem}\label{th_QIZ}
	Let $F$ be the distribution function of the probability law on $\Z^d$. Let $f$ be its characteristic function. Then $F$ is quasi-infinitely divisible if and only if $f(t)\ne 0$ for all $t\in\R^d$. In that case, $f$ admits the following representation
	\begin{eqnarray}\label{th_QIZ_repr}
		f(t)=\exp\biggl\{i\lng t,\gamma\rng+  \sum_{k\in\Z^d\setminus\{\oo\}}\lambda_{k} \bigl(e^{i\lng t,k\rng}-1\bigr)\biggr\},\quad t\in\R^d,
	\end{eqnarray}
	where $\gamma\in\Z^d$, $\lambda_k\in\R$, $k\in\Z^d\setminus\{\oo\}$, and $\sum_{k\in\Z^d\setminus\{\oo\}}|\lambda_{k}|<\infty$. 
\end{Theorem}
It is clear that \eqref{th_QIZ_repr} can be rewritten in the form \eqref{ID_Levy_representation}. Using this theorem, the authors of \cite{Lindner_Berger} also proved the Cram\'er--Wold device for infinite divisibility of $\Z^d$-valued distributions. We formulate the corresponding result in a simplified form omitting equivalent propositions.  
\begin{Theorem}\label{th:Cr-W}
	Let $\xi$ be a random vector on $\Z^d$ with distribution function $F$. Then $F$ is infinitely divisible if and only if for any $c \in \R^d$ distribution functions $F_c$ of random variables $\lng c, \xi\rng$ are infinitely divisible. 
\end{Theorem}

Recall that the classical Cram\'er--Wold device is a fact that a probability distribution of a $d$-dimensional random vector $\xi$ is uniquely determined by distributions of all  linear combinations of its components, i.e. by distributions of $\lng c, \xi \rng$ for all $c \in \R^d$ (see \cite{Cramer_Wold}). Statements concerning some property for multivariate random vectors that can be expressed by corresponding statements for its linear combinations are also called as Cram\'er--Wold devices. So Theorem \ref{th:Cr-W} is an interesting particular example of this.  The Cram\'er--Wold device is well known for strict and symmetric stability (see \cite{Samor_Taqqu}): $d$-dimensional random vector $\xi$ is strictly (symmetric) stable if and only if random variable $\lng c, \xi \rng$ is strictly (symmetric) stable for any $c \in \R^d$. Note that, however, the Cram\'er--Wold device for infinite divisibility does not hold in general. If a $d$-dimensional random vector $\xi$ has infinitely divisible distribution, then the distribution of $\lng c, \xi \rng$ is infinitely divisible, too, for all $c \in \R^d$, but for $d \geqslant 2$ there exist examples that the converse is not true (see \cite{Dwass} and \cite{Ibragimov}). 

The purpose of this article is to generalize Theorems \ref{th_QIZ} and  \ref{th:Cr-W} to arbitrary multivariate discrete distribution functions. 
More precisely, we obtain a criterion of quasi-infinitely divisibility with  similar representation \eqref{th_QIZ_repr} and also we prove the Cram\'er--Wold devices for infinite and quasi-infinite divisibility for such laws. The corresponding results will be formulated in Section 2. The necessary tools, which are also of independent interest,  will be formulated in  Section 3. All of the mentioned results will be proved in Section 4.

\section{Main results}
Let us consider a  multivariate discrete probability law with the following distribution function
\begin{eqnarray}\label{main_rep}
	F(x)= \sum_{\substack{k\in \N:\\ x_k\in(-\infty, x]}} p_{x_k}, \quad x\in \R^d.
\end{eqnarray}
Here  $x_k \in \R^d$, $k \in \N$, are distinct numbers with probability weights $p_{x_k} \geqslant 0$, $k \in \N$ (the set of positive integers), $\sum_{k=1}^{\infty}p_{x_k} = 1$. We denote by $(-\infty, x]$ with $x=(x^{(1)},x^{(2)},\ldots, x^{(d)})\in \R^d$  the set $(-\infty, x^{(1)}]\times\ldots\times   (-\infty, x^{(d)}]\subset \R^d$. Let $f$ be the characteristic function of $F$, i.e.
\begin{eqnarray}\label{ch_fun}
	f(t) \defeq \int_{\R^d} e^{i\lng t,x\rng} \dd F(x)= \sum_{k\in \N} p_{x_k}e^{i\lng t,x_k\rng}, \quad t \in \R^d.
\end{eqnarray}

We will formulate a criterion for the distribution function $F$  to be quasi-infinitely divisible through condition for characteristic function  $f$. For the sharp formulation of the result below we need to introduce the set of all finite $\Z^d$-linear combinations of elements from a set $Y \subset \CC^d$ ($\CC$ is the set of complex numbers):
\begin{eqnarray}\label{modulus}
	\langle Y \rangle = \Bigr\{ \sum_{k=1}^n z_k y^{(1)}_k \colon \, n \in \N, \, z_k \in \Z, \, y_k \in Y\Bigr\} \times \ldots \times\Bigr\{ \sum_{k=1}^n z_k y_k^{(d)} \colon \, n \in \N, \, z_k \in \Z, \, y_k \in Y\Bigr\}, 
\end{eqnarray}
where $y_k = \bigr(y^{(1)}_k, \ldots, y^{(d)}_k\bigr)$.	So $\langle Y\rangle$ is a module over the ring $\Z^d$ with the generating set $Y$. 
It is easily seen that $Y \subset \langle Y \rangle$, $\oo \in \langle Y \rangle$. If a countable set $Y \ne \varnothing$, then $\langle Y \rangle$ is an infinite countable set.

\begin{Theorem}\label{th_Repr}
	Let $F$ be a discrete distribution function of the form \eqref{main_rep} with characteristic function $f$ of the form \eqref{ch_fun}. The following statements are equivalent:
	\begin{enumerate}
		\item[$(a)$] $F$ is quasi-infinitely divisible;
		\item[$(b)$] $\inf_{t\in\R^d} |f(t)| > 0$.
	\end{enumerate}
	If one of the conditions is satisfied, and hence all, then $f$ admits the following representation
	\begin{eqnarray}\label{discr_Levy}
		f(t) = \exp\biggr\{ i\lng t,\gamma \rng + \sum_{u \in \langle X \rangle \setminus\{\oo\}} \lambda_u \bigr(e^{i\lng t,u\rng} - 1\bigr) \biggr\}, \quad t\in \R^d,
	\end{eqnarray} 
	where $X \defeq \{ x_k \colon  p_{x_k} > 0, k \in \N\} \ne \varnothing$, $\gamma \in \lng X \rng $, $\lambda_u \in \R$ for all $u \in \langle X \rangle \setminus\{\oo\}$, and $\sum_{u \in \langle X \rangle \setminus\{0\}}|\lambda_u| < \infty$.
\end{Theorem}

It is easily seen that \eqref{discr_Levy} can be rewritten in the form \eqref{ID_Levy_representation}. So if characteristic function of multivariate discrete probability law is represented by \eqref{discr_Levy}, then the corresponding distribution function $F$ is quasi-infinitely divisible. Also observe that, by this theorem and on account of the conditions and uniqueness of the L\'evy  representation \eqref{ID_Levy_representation}, the multivariate discrete distribution function $F$  is infinitely divisible if and only if its characteristic function $f$ admits representation \eqref{discr_Levy} with the same $X$ and $\gamma$, but with $\lambda_u \geqslant 0$ for all $u \in \langle X \rangle \setminus\{\oo\}$, and $\sum_{u \in \langle X \rangle \setminus\{0\}}\lambda_u < \infty$.  

Note that Theorem~\ref{th_Repr}  generalizes Theorem \ref{th_QIZ}. Indeed, for characteristic function $f$ of probability law  on $\Z^d$ the condition that $f(t)\ne 0$, $t\in\R^d$,  is equivalent to the condition that $\inf_{t\in\R}|f(t)|>0$. It follows due to the continuity and $2\pi$-periodicity of the function $|f(t)|$, $t=(t_1,\ldots,t_d)\in\R^d$, over each $t_j$. Theorem~\ref{th_Repr} also  generalizes  the corresponding results from  \cite{Alexeev_Khartov} and \cite{Khartov2} for the discrete distributions  in the univariate case.

We now formulate the Cram\'er--Wold devices for the infinite and quasi-infinite divisibility of multivariate discrete distribution functions. 

\begin{Theorem}\label{th:Cramer-Wold_QID}
Let $\xi$ be a discrete random vector with distribution function $F$ of the form \eqref{main_rep}. Then $F$ is $($quasi-$)$infinitely divisible if and only if for any $c \in \R^d$ the distribution functions $F_c$ of random variables $\lng c, \xi \rng$ are $($quasi-$)$infinitely divisible. 
\end{Theorem}

It is seen that Theorem \ref{th:Cr-W} directly follows from the Theorem \ref{th:Cramer-Wold_QID}. It should be noted that Theorem \ref{th:Cramer-Wold_QID} does not contradict with the results from the \cite{Dwass} and \cite{Ibragimov}, because the distributions from the counterexamples contained the absolutely continuous part. 

\section{Tools}

We will get the main result from more general positions. Namely, we will consequently study admission of the L\'evy type representations for a general almost periodic function $h$ on $\R^d$ that is very similar to  $f$.  

\begin{Theorem}\label{suf_cond_discr}
	Let  $h\colon \R^d\to \CC$ be a function of the following form:
	\begin{eqnarray*}
		h(t) = \sum_{y \in Y} q_y e^{i\lng t,y\rng},\quad t\in \R^d,
	\end{eqnarray*}
	where $Y\subset \R^d$ is a nonempty at most countable set, $q_y \in \CC$ for all $y\in Y$, and $0<\sum_{y \in Y}|q_y|<\infty$. Assume that $h(\oo) =\sum_{y \in Y}q_y =1$. If $\inf_{t\in \R^d}|h(t)| = \mu > 0$, then $h$ admits the following representation
	\begin{eqnarray}\label{alm_per_Levy}
		h(t) = \exp\biggr\{ i\lng t,\gamma\rng + \sum_{u \in \lng Y \rng \setminus\{\oo\}} \lambda_u \bigr(e^{i\lng t,u\rng} - 1\bigr) \biggr\}, \quad t\in \R^d,
	\end{eqnarray} 
	where $\gamma \in \langle Y \rangle$, $\lambda_u \in \CC$ for all $u\in \langle Y \rangle \setminus\{\oo\}$, and $\sum_{u \in \lng Y \rng \setminus\{\oo\}}|\lambda_u| < \infty$.
\end{Theorem}

It should be noted that the function $h$ in Theorem \ref{suf_cond_discr} is an \textit{almost periodic function} on $\R^n$ with absolutely convergence Fourier series. Recal that (see \cite[p. 255]{Levitan} or \cite[Definition 1]{Neumann}) a function $h\colon \R^d \to \CC$ is called almost periodic if for any sequence $\{t_n\}_{n\in \N}$ from $\R^d$ there exists a subsequence $(t_{n_k})_{k\in\N}$ and a continuous function $\varphi\colon \R^d \to \CC$ such that 
\begin{eqnarray*}
	\sup_{t\in \R^d}\bigr|h(t+t_{n_k}) - \varphi(t)\bigr| \underset{k \to \infty} \longrightarrow 0. 
\end{eqnarray*}
Detailed information about almost periodic functions on $\R^n$ can be found in \cite{BalanKrishtal}, \cite{Levitan}, \cite{LevitanZhikov}, \cite{Neumann}, and \cite{Pankov}. Note that in the above literature, the results are formulated in a greater generality, that is, for local compact Abelian (LCA) groups.

We now turn to the following general version of the Theorem \ref{th_Repr}. 
\begin{Theorem}\label{alper_mainres}
	Let  $h\colon \R^d\to \CC$ be a function of the following form
	\begin{eqnarray*}
		h(t)=\sum_{y\in Y} q_y e^{i\lng t,y\rng},\quad t\in\R^d,
	\end{eqnarray*}
	where $Y\subset \R^d$ is a nonempty at most countable set, $q_y \in \CC$ for all $y\in Y$, and $0<\sum_{y \in Y}|q_y|<\infty$. Suppose that $h(\oo) =\sum_{y\in Y} q_{y} =1$. Then the following statements are equivalent:
	\begin{enumerate}
		\item[$(i)$] $\inf_{t \in \R^d}|h(t)|  >0;$
		\item[$(ii)$] There exist a countable set $Z \subset \R^d$ and coefficients $r_{z} \in \CC$, $z\in Z$, $\sum_{z \in Z}|r_{z}|<\infty$, such that
		\begin{eqnarray*}
			\frac{1}{h(t)} = \sum_{z \in Z}r_z e^{i\lng t,z\rng},\quad t\in\R^d;
		\end{eqnarray*}
		\item[$(iii)$] $h$ admits the representation
		\begin{eqnarray*}
			h(t) = \exp\biggl\{ i\lng t,\gamma\rng + \sum_{u \in \langle Y \rangle \setminus\{\oo\}} \lambda_u \bigr(e^{i\lng t,u\rng} - 1\bigr) \biggr\}, \quad t\in \R^d,
		\end{eqnarray*} 
		where $\gamma \in \langle Y \rangle$, $\lambda_u \in \CC$ for all $u \in \langle Y \rangle \setminus\{\oo\}$, and $\sum_{u \in \langle Y \rangle \setminus\{0\}}|\lambda_u| < \infty;$ 
		\item[$(iv)$] $h$ admits the representation
		\begin{eqnarray}\label{QID_fol_general_case}
			h(t) = \exp\biggr\{ i\lng t,\gamma\rng - \tfrac{1}{2}\lng t, Qt\rng + \int_{\R^d} \Bigr(e^{i\lng t,u\rng} - 1 - \tfrac{i\lng t,u\rng}{1+\|u\|^2}\Bigr) \nu(\dd u) \biggr\}, \quad t\in \R^d,
		\end{eqnarray} 
		where $\gamma \in \CC^d$,  $Q \in \CC^{d \times d}$ is a matrix, $\nu$ is a complex measure on $\R^d$ such that  $$\nu(\{\oo\}) = 0,\quad\text{and}\quad\int_{\R^d}\min\bigl\{\|x\|^2, 1\bigr\}|\nu|(\dd x)<\infty.$$
	\end{enumerate}
\end{Theorem}

\section{Proofs}

\textbf{Proof of Theorem \ref{suf_cond_discr}.}  We will sequentially consider the following cases: 1)  $Y = \Z^d$, 2)   $Y$ is a  finite subset of $\R^d$,  3) $Y$ is at  most countable subset of $\R^d$ (the general case). We  always assume that  $Y\ne \varnothing$. Each subsequent case will be based on the previous one.

\textbf{1)} Suppose that $Y = \Z^d$. It is easy to see that the function $h$ is $2\pi$-periodic in all coordinates, i.e. for any $k=1, \ldots, d$ and  $t \in \R^d$ we have $h(t + 2\pi e_k) = h(t)$, where $\{e_1, e_2,\ldots, e_d \}$ denotes the canonical basis in $\R^d$. 
Let us consider the distinguished logarithm $t\mapsto\Ln h(t)$,  $t\in\R^d$, which satisfies $\exp\{ \Ln h(t)\}=h(t)$, $t\in\R^d$, and it is uniquely defined by continuity with the condition  $\Ln h(\oo)=0$ (see \cite[Lemma 7.6]{Sato_book}). For any $k = 1, \ldots, d$ we have
\begin{eqnarray*}
	\exp\{ \Ln h(t+2\pi e_k)\} = h(t+2\pi e_k) = h(t) = \exp\{ \Ln h(t)\},\quad t\in\R^d. 
\end{eqnarray*}
So $\Ln h(t+2\pi e_k) - \Ln h(t) \in 2\pi i \Z$ for any $k = 1, \ldots, d$ and $t \in \R^d$. Since $t \mapsto \Ln h(t+2\pi e_k) - \Ln h(t)$ is a continuous function on $\R^d$, there exist constants $\gamma_1,\ldots, \gamma_d \in \Z$ such that
\begin{eqnarray*}
	\gamma_k = \frac{\Ln h(t+2\pi e_k) - \Ln h(t)}{2\pi i},\quad t\in\R^d, \quad k = 1, \ldots, d.
\end{eqnarray*} 
Let us define the vector $\gamma = (\gamma_1, \ldots, \gamma_d)^T\in \Z^d$. So the function $t\mapsto \Ln h(t) - i\lng t, \gamma\rng$ is $2\pi$-periodic in all coordinates. By  \cite[Proposition 3.1]{Lindner_Berger}, one can conclude that
\begin{eqnarray*}
	\Ln h(t) = i\lng t, \gamma\rng + \sum_{u \in \Z^d \setminus\{\oo\}} \lambda_u \bigr(e^{i\lng t, u\rng} - 1\bigr),\quad t\in\R^d,
\end{eqnarray*} 
where $\lambda_u \in \CC$ for all $u \in \Z^d \setminus\{\oo\}$, and  $\sum_{u \in \Z^d \setminus\{\oo\}} |\lambda_u| < \infty$. Note that in this case $\langle Y \rangle = \Z^d$.

\textbf{2)} Assume  that $Y=\{y_1, \ldots, y_n\}$, where $y_1, \ldots, y_n$ are distinct elements from $\R^d$.  So we have $h(t) = \sum_{k=1}^n q_{y_k} e^{i\lng t,y_k\rng}$, $t\in \R^d$. If $n=1$ then $Y=\{y_1\}$ and $q_{y_1}=1$. For this case representation \eqref{alm_per_Levy} holds with $\gamma=y_1$   and $\lambda_{u}=0$ for all $u\in \langle Y \rangle \setminus\{\oo\}$. We next suppose that $n\geqslant 2$.  We set $y_k=(y_k^{(1)},\ldots, y_k^{(d)})$, $k=1,\ldots, n$. Without loss of generality, we can assume that for every $j=1,\ldots, d$ there exist  $k=1,\ldots, n$ such that $y^{(j)}_k\ne 0$, since otherwise we can turn to the space $\R^{d'}$ with some $d'<d$. Next, for every $j = 1, \ldots, d$ we can choose non-zero $\beta^{(j)}_1, \ldots, \beta^{(j)}_{m_j} \in  Y^{(j)} = \{ y_1^{(j)}, \ldots, y_n^{(j)}\} \subset \R$ that consist a basis in $Y^{(j)}$ over $\Z$, i.e. for every $j = 1, \ldots, d$ and for every $k = 1, \ldots, n$ there exist unique values $c^{(j)}_{k,1},\ldots, c^{(j)}_{k,{m_j}} \in \Z$ such that $y^{(j)}_k = \sum_{l= 1}^{m_j} c^{(j)}_{k,l} \beta^{(j)}_l$. So it is easy to check that 
\begin{eqnarray*}
	\langle Y \rangle = \Bigr\{ \sum_{l=1}^{m_1} z_l^{(1)} \beta^{(1)}_l \colon \, z_l^{(1)} \in \Z\Bigr\} \times \ldots \times\Bigr\{ \sum_{l=1}^{m_d} z_l^{(d)} \beta^{(d)}_l \colon \, z_l^{(d)} \in \Z\Bigr\}.  
\end{eqnarray*}
Note that for every $j = 1, \ldots, d$ the values $\beta^{(j)}_1, \ldots, \beta^{(j)}_{m_j}$ are linearly independent over $\Z$, i.e. the equation $l_1 \beta^{(j)}_1 + \ldots+ l_{m_j}\beta^{(j)}_{m_j} = 0$ holds with $l_1, \ldots, l_{m_j} \in \Z$  if and only if $l_1 = \ldots= l_{m_j} = 0$. 	

We now consider the function
\begin{eqnarray}\label{def_phi}
	\varphi\bigr(t^{(1)}_1, \ldots, t^{(1)}_{m_1}, \ldots, t^{(d)}_1, \ldots, t^{(d)}_{m_d}\bigr) \defeq \sum_{k=1}^n q_{y_k} \exp\biggr\{i\sum_{j = 1}^d\sum_{l=1}^{m_j} c^{(j)}_{k,l} \beta^{(j)}_l t^{(j)}_l\biggr\}, 
\end{eqnarray}
where $t^{(j)}_l \in \R$, $l = 1, \ldots, m_j$, and $j = 1, \ldots, d$. If for any such $j$ and $l$ we set $t^{(j)}_l \defeq t^{(j)}\in \R$, then  
\begin{eqnarray}\label{diagonal}
	\varphi\bigr(t^{(1)}_1, \ldots, t^{(1)}_{m_1}, \ldots, t^{(d)}_1, \ldots, t^{(d)}_{m_d}\bigr) = h(t), 
\end{eqnarray}
where $t=(t^{(1)}, \ldots, t^{(d)})$.	We set $M \defeq m_1 + \ldots +m_d$. Let us fix an arbitrary $\e >0$. Since the function $\varphi$  is uniformly continuous, there exists $\delta_\e > 0$ such that for any $\tilde t_1$ and $\tilde t_2$ from $\R^M$ satisfying $\| \tilde t_1 - \tilde t_2\| < \delta_e$  we have $\bigr|\varphi(\tilde t_1) - \varphi(\tilde t_2)\bigr| < \e$. Let us arbitrarily fix  the vector $t\defeq \bigl(t^{(1)}_1, \ldots, t^{(1)}_{m_1}, \ldots, t^{(d)}_1, \ldots, t^{(d)}_{m_d}\bigr)\in\R^M$. We set $b_j\defeq \min \bigl\{|\beta^{(j)}_1|,\ldots, |\beta^{(j)}_{m_j}|\bigr\}>0$ for every $j = 1, \ldots, d$. Since for every $j$ the values $\beta_1^{(j)}, \ldots, \beta_{m_j}^{(j)}$ are linearly independent over $\Z$, then, by the Kronecker theorem (see \cite[p.37]{LevitanZhikov}), we conclude that the inequalities
\begin{eqnarray*}
	\bigr| \beta_l^{(j)} s^{(j)} - t_l^{(j)}  - 2\pi n_l^{(j)}\bigr| < \tfrac{\delta_\e b_j}{\sqrt{m_j d}},\quad l = 1, \ldots, m_j,
\end{eqnarray*}
have a common solution $s^{(j)} \in \R$  for some $n_l^{(j)} \in \Z$. We fix these numbers  and we conclude that
\begin{eqnarray*}
	\Biggr|  s^{(j)} - \dfrac{t_l^{(j)} + 2\pi n_l^{(j)}}{\beta_l^{(j)}}\Biggr| < \dfrac{\delta_\e}{\sqrt{m_j d}},\quad l = 1, \ldots, m_j,
\end{eqnarray*}
and
\begin{eqnarray*}
	\sum_{j = 1}^d\sum_{l=1}^{m_j}\Biggr|  s^{(j)} - \dfrac{t_l^{(j)} + 2\pi n_l^{(j)}}{\beta_l^{(j)}}\Biggr|^2 < \delta_\e^2.
\end{eqnarray*}
The latter inequality means that $\|s-\tilde t\, \|<\delta_\e$, where   
\begin{eqnarray*}
	&&s\defeq \bigl(s^{(1)},\ldots, s^{(1)}, \ldots, s^{(d)},\ldots, s^{(d)} \bigr)\in \R^M,\\
	&&\tilde t\defeq\biggl(\tfrac{t^{(1)}_1 + 2\pi n_1^{(1)}}{\beta_1^{(1)}}, \ldots, \tfrac{t^{(1)}_{m_1}+ 2\pi n_{m_1}^{(1)}}{\beta_{m_1}^{(1)}}, \ldots, \tfrac{t^{(d)}_1+ 2\pi n_1^{(d)}}{\beta_1^{(d)}}, \ldots, \tfrac{t^{(d)}_{m_d}+ 2\pi n_{m_d}^{(d)}}{\beta_{m_d}^{(d)}}\biggr)\in \R^M, 
\end{eqnarray*}
in the vector $s$: $s^{(1)}$ repeats $m_1$ times, $s^{(2)}$ repeats $m_2$ times, \ldots, $s^{(d)}$ repeats $m_d$ times. Therefore $\bigr|\varphi(s) - \varphi(\tilde t)\bigr| < \e$. 	It is easily seen  from \eqref{def_phi} that
\begin{eqnarray*}
	\varphi(\tilde t )&=&\varphi\biggl(\tfrac{t^{(1)}_1 + 2\pi n_1^{(1)}}{\beta_1^{(1)}}, \ldots, \tfrac{t^{(1)}_{m_1}+ 2\pi n_{m_1}^{(1)}}{\beta_{m_1}^{(1)}}, \ldots, \tfrac{t^{(d)}_1+ 2\pi n_1^{(d)}}{\beta_1^{(d)}}, \ldots, \tfrac{t^{(d)}_{m_d}+ 2\pi n_{m_d}^{(d)}}{\beta_{m_d}^{(d)}}\biggr)\\
	&=& \varphi\biggl(\tfrac{t^{(1)}_1}{\beta_1^{(1)}}, \ldots, \tfrac{t^{(1)}_{m_1}}{\beta_{m_1}^{(1)}}, \ldots, \tfrac{t^{(d)}_1}{\beta_1^{(d)}}, \ldots, \tfrac{t^{(d)}_{m_d}}{\beta_{m_d}^{(d)}}\biggr)\\
	&=:&\tilde{\varphi} (t),
\end{eqnarray*}
i.e.
\begin{eqnarray}\label{def_tildephi}
	\tilde{\varphi} (t)=\sum_{k=1}^n q_{y_k} \exp\biggr\{i\sum_{j = 1}^d\sum_{l=1}^{m_j} c^{(j)}_{k,l} t^{(j)}_l\biggr\};
\end{eqnarray}
since  $t$ was fixed arbitrarily, we consider $\tilde{\varphi}$ as a function on $\R^M$.	So we have that $\bigr|\varphi(s) - \tilde\varphi(t)\bigr| < \e$. Thus,  due to \eqref{diagonal}, we get that for any $\e>0$ and  $t\in \R^M$ there exists $ s'=\bigl(s^{(1)},\ldots, s^{(d)} \bigr) \in \R^d$ such that  $\bigr|h(s') - \tilde\varphi(t)\bigr| < \e$. According to the assumption  $\inf_{s \in \R^d}|h(s)| > 0$, we conclude that $\inf_{t \in \R^M}|\tilde\varphi(t)| > 0$. 

We now apply the previous part \textbf{1)}  to the function \eqref{def_tildephi} (it is valid, because there are $c^{(j)}_{k,l}\in\Z$ in \eqref{def_tildephi}). So we have the following representation: 
\begin{eqnarray*}
	\Ln\tilde\varphi(t)&=& \Ln \tilde{\varphi}\bigr(t^{(1)}_1, \ldots, t^{(1)}_{m_1}, \ldots, t^{(d)}_1, \ldots, t^{(d)}_{m_d}\bigr)\\
	&=& i\sum_{j=1}^d\sum_{l=1}^{m_j} \gamma_l^{(j)} t_l^{(j)} + \sum_{ z \in \Z^M\setminus\{\oo\}} \lambda_{z} \biggl( \exp\biggl\{i\sum_{j=1}^d\sum_{l=1}^{m_j} z_l^{(j)} t_l^{(j)} \biggr\} - 1\biggr),
\end{eqnarray*}
where  $z = \bigr(z^{(1)}_1, \ldots, z^{(1)}_{m_1}, \ldots, z^{(d)}_1, \ldots, z^{(d)}_{m_d}\bigr) \in \Z^M\setminus\{\oo\}$, $\gamma_l^{(j)} \in \Z$, $\lambda_{z} \in \CC$ for all $z \in \Z^M\setminus\{\oo\}$, and $\sum_{z \in \Z^M\setminus\{\oo\}}|\lambda_{z}| < \infty$. From the above, we get
\begin{eqnarray*}
	\Ln \varphi\bigr(t^{(1)}_1, \ldots, t^{(1)}_{m_1}, \ldots, t^{(d)}_1, \ldots, t^{(d)}_{m_d}\bigr) &=& i\sum_{j=1}^d\sum_{l=1}^{m_j} \gamma_l^{(j)}\beta_l^{(j)} t_l^{(j)} \\&+& \sum_{ z \in \Z^M\setminus\{\oo\}} \lambda_{z} \Biggl( \exp\biggl\{i\sum_{j=1}^d\sum_{l=1}^{m_j} z_l^{(j)}\beta_l^{(j)} t_l^{(j)} \biggr\} - 1\Biggr).
\end{eqnarray*}
Due to \eqref{diagonal}, for every $t = (t^{(1)}, \ldots, t^{(d)})$ we have
\begin{eqnarray*}
	\Ln h(t) = i\sum_{j=1}^d\biggl(\sum_{l=1}^{m_j} \gamma_l^{(j)}\beta_l^{(j)}\biggr) t^{(j)} + \sum_{z \in \Z^M\setminus\{\oo\}} \lambda_{z} \Biggl( \exp\biggl\{i\sum_{j=1}^d\biggl(\sum_{l=1}^{m_j} z_l^{(j)}\beta_l^{(j)}\biggr) t^{(j)} \biggr\} - 1\Biggr).
\end{eqnarray*}
For every $j=1, \ldots, d$ we set $\gamma^{(j)} \defeq \sum_{l=1}^{m_j} \gamma_l^{(j)}\beta_l^{(j)} \in \langle Y^{(j)}\rangle$, $\gamma \defeq  (\gamma^{(1)}, \ldots, \gamma^{(d)}) \in \lng Y\rng$, $\lambda_u \defeq  \lambda_{z}$ for $u = (u^{(1)}, \ldots, u^{(d)}) \in \lng Y\rng\setminus\{\oo\}$ with  $u^{(j)} \defeq \sum_{l=1}^{m_j} z_l^{(j)}\beta_l^{(j)}\in \langle Y^{(j)} \rangle \setminus\{0\}$ ($u$ determines $z$ uniquely, because $\beta_l^{(j)}$ consist a basis, see above). Thus we come  to the representation \eqref{alm_per_Levy} for $h$. 

\textbf{3)} We now turn to the general case: $Y$ is at most countable subset of $\R^d$. Without loss of generality we can set $Y\defeq\{y_1,y_2,\ldots\}$ with distinct $y_k\in\R^d$. So $A\defeq\sum_{k=1}^{\infty} |q_{y_k}|<\infty$ and $h(t) \defeq \sum_{k=1}^{\infty} q_{y_k} e^{i\lng t,y_k\rng}$, $t \in \R^d$. 
We will  approximate $h$ by the following functions:
\begin{eqnarray*}
	h_n(t)\defeq \sum_{k=1}^{n}  q_{n,y_k} e^{i\lng t,y_k\rng},\quad  t\in\R^d,\quad n\in\N,
\end{eqnarray*}
where
\begin{eqnarray*}
	q_{n,y_k}\defeq \dfrac{q_{y_k}}{\sum_{m=1}^{n} q_{y_m}},\quad k=1,\ldots, n,\quad n\in\N.
\end{eqnarray*}
Since $\sum_{k=1}^{\infty} q_{y_k}=1$, here $\bigl|\sum_{m=1}^{n} q_{y_m}\bigr|\geqslant \tfrac{1}{2}$ for all $n\geqslant n_0$ with a positive integer $n_0$.  Let us estimate the approximation error for every $n\geqslant n_0$:
\begin{eqnarray*}
	\sup\limits_{t\in\R^d}|h(t)-h_n(t)|&=&\sup\limits_{t\in\R^d}\biggl| \sum_{k=1}^n (q_{y_k}-q_{n,y_k})e^{i\lng t,y_k\rng} +\sum_{k=n+1}^\infty q_{y_k}e^{i\lng t,y_k\rng}\biggr|\\
	&\leqslant& \sum_{k=1}^n |q_{y_k}-q_{n,y_k}|+\sum_{k=n+1}^\infty |q_{y_k}|.
\end{eqnarray*}
Due to $\sum_{m=1}^{\infty} q_{y_m}=1$, we have
\begin{eqnarray*}
	\sum_{k=1}^n |q_{y_k}-q_{n,y_k}|= \biggl|1-\dfrac{1}{\sum_{m=1}^{n} q_{y_m}}\biggr|\cdot\sum_{k=1}^n |q_{y_k}|
	=\biggl|\dfrac{\sum_{m=n+1}^{\infty} q_{y_m}}{\sum_{m=1}^{n} q_{y_m}}\biggr|\cdot \sum_{k=1}^n |q_{y_k}| \leqslant 2A\sum_{m=n+1}^{\infty} |q_{y_m}|.
\end{eqnarray*}
We used  $\sum_{k=1}^{\infty} |q_{y_k}|=A$ and $\bigl|\sum_{m=1}^{n} q_{y_m}\bigr|\geqslant \tfrac{1}{2}$ in the last inequality. Thus we obtain
\begin{eqnarray*}
	\sup\limits_{t\in\R^d}|h(t)-h_n(t)|\leqslant (2A+1)\sum_{m=n+1}^\infty |q_{y_m}|,\quad n\geqslant n_0.
\end{eqnarray*}
Since $\sum_{k=1}^n |q_{y_k}| < \infty$, we have that $\sup_{t \in \R^d}|h(t)-h_n(t)|\to 0$, $n\to\infty$. Hence for any fixed $\e \in (0,\tfrac{1}{4})$ there exists a positive integer $n_{\e} \geqslant n_0$ such that for every $n \geqslant n_{\e}$ we have
\begin{eqnarray}\label{upper_esthhn}
	\sup_{t \in \R^d} \bigr| h(t) - h_n(t) \bigr| \leqslant \e\mu,
\end{eqnarray}
where we set $\mu \defeq \inf_{t\in \R^d}|h(t)| >0$. So for every $n \geqslant n_{\e}$
\begin{eqnarray}\label{lower_esthn}
	\inf_{t \in \R^d} \bigr| h_n(t) \bigr| \geqslant \inf_{t \in \R^d} \bigr| h(t)\bigr| - \sup_{t \in \R^d} \bigr| h(t) - h_n(t) \bigr|  \geqslant (1-\e)\mu. 
\end{eqnarray}

We now fix $n\geqslant n_\e$ and we represent $h(t)=h_n(t)\cdot R_n(t)$ with   $R_n(t)\defeq h(t)/h_n(t)$, $t\in\R^d$. Since $h$, $h_n$, $R_n$ are continuous functions without zeroes on $\R^d$  and they  equal $1$ at $t=\oo$, we can proceed to the distinguished logarithms: 
\begin{eqnarray}\label{eq_Lnh}
	\Ln h(t)=\Ln h_n(t) + \Ln R_n(t),\quad t\in\R^d.
\end{eqnarray}

Let us consider the function $\Ln h_n$. By the result of part \textbf{2)}, we have
\begin{eqnarray*}
	\Ln h_n(t)=i\lng t,\gamma_n\rng+
	\sum_{u\in\langle Y_n\rangle\setminus\{\oo\}}\lambda_{n,u} \bigl(e^{i\lng t,u\rng}-1\bigr),\quad t\in\R^d,
\end{eqnarray*}
with a set $Y_n\defeq\bigl\{y_k\colon\, q_{y_k}\ne0,\, k=1,\ldots,n \bigr\}$, and numbers $\gamma_n\in\langle Y_n\rangle$, $\lambda_{n,u}\in\CC$ for all $u\in\langle Y_n \rangle\setminus\{\oo\}$, $\sum_{u\in\langle Y_n\rangle\setminus\{\oo\}} |\lambda_{n,u}|<~\infty$. Setting $\lambda_{n,\oo}\defeq -\sum_{u\in\langle Y_n\rangle\setminus\{\oo\}} \lambda_{n,u} \in \CC$ , we represent $\Ln h_n$ in the following form
\begin{eqnarray*}
	\Ln h_n(t)=i\lng t,\gamma_n\rng+
	\sum_{u\in\langle Y_n\rangle}\lambda_{n,u} e^{i\lng t,u\rng},\quad t\in\R^d.
\end{eqnarray*}
Observe that $Y_n\subset Y$, and hence $\langle Y_n\rangle\subset \langle Y\rangle$. So we can write
\begin{eqnarray}\label{eq_fn_moduleX}
	\Ln h_n(t)=i\lng t,\gamma_n\rng+
	\sum_{u\in\langle Y\rangle}\lambda_{n,u} e^{i\lng t,u\rng},\quad t\in\R^d,
\end{eqnarray}
where for every $u\in \langle Y\rangle\setminus \langle Y_n\rangle$ we define $\lambda_{n,u}\defeq 0$  for the case $\langle  Y\rangle\setminus \langle Y_n\rangle\ne \varnothing$. 

We next consider the function $\Ln R_n$. Observe that
\begin{eqnarray}\label{eq_Lnhhn}
	\Ln R_n(t)=\ln\biggl(1+\dfrac{h(t)-h_n(t)}{h_n(t)}\biggr),\quad t\in\R^d,
\end{eqnarray}
where the latter is the principal value of the logarithm. Indeed,  due to \eqref{upper_esthhn} and \eqref{lower_esthn}, 
\begin{eqnarray}\label{conc_supt_hhn}
	\sup_{t\in\R^d} \biggl|\dfrac{h(t) -h_n(t)}{h_n(t)}\biggr| \leqslant \dfrac{\sup_{t\in\R^d}|h(t) -h_n(t)|}{\inf_{t\in\R^d}|h_n(t)|}\leqslant \dfrac{\e}{1-\e}<1,
\end{eqnarray}
and the function in the right-hand side of \eqref{eq_Lnhhn} is continuous and it equals $0$ at $t=\oo$.
Therefore we get the decomposition
\begin{eqnarray*}
	\Ln R_n(t)=\sum\limits_{m=1}^{\infty} \dfrac{(-1)^{m-1}}{m} \biggl(\dfrac{h(t)-h_n(t)}{h_n(t)}\biggr)^m,\quad t\in\R^d,
\end{eqnarray*}
which yields the estimate
\begin{eqnarray*}
	\sup_{t\in\R}|\Ln R_n(t)|\leqslant \sum\limits_{m=1}^{\infty} \dfrac{1}{m}\, \sup_{t\in\R^d}\biggl|\dfrac{h(t)-f_n(t)}{f_n(t)}\biggr|^m\leqslant \sum\limits_{m=1}^{\infty} \dfrac{1}{m} \biggl(\dfrac{\e}{1-\e}\biggr)^m.
\end{eqnarray*}
Since $\e\in(0,\tfrac{1}{4})$, we have
\begin{eqnarray*}
	\sum\limits_{m=1}^{\infty} \dfrac{1}{m} \biggl(\dfrac{\e}{1-\e}\biggr)^m\leqslant\sum\limits_{m=1}^{\infty} \biggl(\dfrac{\e}{1-\e}\biggr)^m=\dfrac{\tfrac{\e}{1-\e}}{1-\tfrac{\e}{1-\e}}=\dfrac{\e}{1-2\e}<2\e.
\end{eqnarray*}
Thus we obtain
\begin{eqnarray}\label{ineq_supRn}
	\sup_{t\in\R^d}|\Ln R_n(t)|<2\e.
\end{eqnarray}

Let us consider the function $(h-h_n)/h_n$ from \eqref{eq_Lnhhn}. It is clear that $h-h_n$ is an almost periodic function with absolutely convergent Fourier series. Due to \cite[Theorem 3.2]{BalanKrishtal} the function $1/h_n$ is also an almost periodic with absolutely convergent Fourier series. Since the function $z\mapsto\ln(1+z)$, $z\in\CC$, is analytic on the unit disk, due to \eqref{conc_supt_hhn}, \cite[Theorem 3.2]{BalanKrishtal}, and \cite[Corollary 5.15]{Grochenig}, we get that $\Ln R_n$ is an almost periodic function with absolutely convergent Fourier series:
\begin{eqnarray}\label{eq_LnRn}
	\Ln R_n(t)=\sum\limits_{u\in \Delta_n} \beta_{n,u} e^{i\lng t,u\rng},\quad t\in\R^d,
\end{eqnarray}
where $\Delta_n$ is at most countable set of vectors from $\R^d$, $\beta_{n,u}\in \CC$ for $u\in \Delta_n$, and $\sum_{u\in \Delta_n}|\beta_{n,u}|<\infty$. 

We now return to the function $\Ln h$ and \eqref{eq_Lnh}.  The formulas \eqref{eq_fn_moduleX} and \eqref{eq_LnRn} yield 
\begin{eqnarray*}
	\Ln h(t)=i\lng t,\gamma_n\rng +\sum_{u\in\lng Y\rng}\lambda_{n,u} e^{i\lng t,u\rng}+ \sum\limits_{u\in \Delta_n} \beta_{n,u} e^{i\lng t,u\rng},\quad t\in\R^d.
\end{eqnarray*}
This formula is valid for every $n\geqslant n_\e$. Let us fix $\mathbf{e} \in \mathbb{S}^{d-1} = \{x \in \R^d\colon\, \|x\| = 1\}$, that is the unit sphere in $\R^d$. 	Since $\sum_{u\in\lng Y\rng}|\lambda_{n,u}|<\infty$ and $\sum_{u\in \Delta_n} |\beta_{n,u}|<\infty$, $n\geqslant n_\e$, it is easy to see that
\begin{eqnarray*}
	\lim\limits_{T \to \infty} \frac{\Ln h(T\mathbf{e})}{iT} = \lng\gamma_n, \mathbf{e}\rng,\quad n\geqslant n_\e.
\end{eqnarray*}
Since the vector $\mathbf{e}$ is choosen arbitrarily from $\mathbb{S}^{d-1}$, $\gamma_n$ are equal for $n\geqslant n_\e$, and  we set $\gamma\defeq \gamma_n$. Due to $\gamma_n\in\lng Y_n\rng \subset \lng Y\rng$, we have $\gamma \in \lng Y \rng$. Thus for every $n \geqslant n_{\e}$ we obtain
\begin{eqnarray*}
	\Ln h(t) = i\lng t,\gamma\rng + \sum_{u\in\langle Y\rangle}\lambda_{n,u} e^{i\lng t,u\rng}  + \sum_{u \in \Delta_n} \beta_{n,u} e^{i\lng t,u\rng},\quad t\in\R^d.
\end{eqnarray*}
Due to the uniqueness theorem for Fourier coefficients (see \cite[Lemma 3.1]{BalanKrishtal}), one can conclude that 
\begin{eqnarray*}
	\Ln h(t) = i\lng t,\gamma\rng + \sum_{u\in\langle Y\rangle}\lambda_{u} e^{i\lng t,u\rng}  + \sum_{u \in Z} \lambda_{u} e^{i\lng t,u\rng},\quad t\in\R^d,
\end{eqnarray*}
where $Z$ is at most countable subset of $\R^d$ such that $\langle Y \rangle \cap Z = \varnothing$, $\lambda_u \in \CC$ for all $u \in \langle Y \rangle \cup Z$, $\sum_{u \in \langle Y \rangle \cup Z}|\lambda_{u}| < \infty$. So for every $n \geqslant n_{\e}$ the following estimate is true:
\begin{eqnarray*}
	\sum_{u \in Z} |\lambda_u|^2 \leqslant \sum_{u \in \Delta_n} |\beta_{n,u}|^2. 
\end{eqnarray*}
Using the Parseval identity (see \cite[Ch. VI, \S 4]{Levitan} or \cite[Theorem 28]{Neumann}) and \eqref{ineq_supRn}, we  get
\begin{eqnarray*}
	\sum_{u \in Z} |\lambda_u|^2 \leqslant \lim\limits_{T\to\infty} \dfrac{1}{(2T)^d} \int\limits_{[-T, T]^d} \bigr|\Ln R_n(t)\bigr|^2\dd t < (2\e)^2,\quad n\geqslant n_\e.
\end{eqnarray*}
Since $\e>0$ can be chosen arbitrarily small, we conclude that $Z = \varnothing$ or $Z \ne \varnothing$, but  $\lambda_u=0$ for all $u\in Z$. Thus
\begin{eqnarray*}
	\Ln h(t) = i\lng t,\gamma\rng + \sum_{u\in\langle Y\rangle}\lambda_{u} e^{i\lng t,u\rng},\quad t\in\R^d,
\end{eqnarray*}	
with $\gamma\in\lng Y\rng$, $\lambda_u\in\CC$ for all $u\in\lng Y\rng$,  and $\sum_{u\in\lng Y\rng}|\lambda_{u}|<\infty$. According to $\bigl(\Ln h(t)-i\lng t,\gamma\rng\bigr)\bigr|_{t=\oo}=0$,  we get $\lambda_{\oo} =-\sum_{u\in\lng Y\rng\setminus\{\oo\}}\lambda_{u}$ and we come to the required representation \eqref{alm_per_Levy}.\quad $\Box$\\

We now return to the proof of the Theorem \ref{alper_mainres}.\\
\textbf{Proof of Theorem \ref{alper_mainres}.} The proof will be carried out in the following sequence: $(ii) \overset{\textbf{I}} \longrightarrow (i) \overset{\textbf{II}} \longrightarrow (iii) \overset{\textbf{III}} \longrightarrow (iv) \overset{\textbf{IV}} \longrightarrow (i) \overset{\textbf{V}} \longrightarrow (ii)$. 

$\textbf{I}.$ Due to $(ii)$, we have 
\begin{eqnarray*}
	\sup_{t\in \R^d} \left| \dfrac{1}{h(t)}\right| \leqslant \sum_{z \in Z}|r_z| = \frac1{\mu} < \infty.
\end{eqnarray*}
It follows that 
\begin{eqnarray*}
	\inf\limits_{t \in \R^d}|h(t)| = \frac{1}{\sup_{t\in \R^d} \left| 1/h(t)\right|} = \mu > 0.
\end{eqnarray*}

$\textbf{II}.$ This implication directly follows  from Theorem \ref{suf_cond_discr}. 

$\textbf{III}.$ It is clear that   $(iii)$ yields $(iv)$ with zero matrix $Q$ and the signed measure
\begin{eqnarray*}
	\nu(B) = \sum_{u \in B\cap \lng Y \rng \setminus\{\oo\}}\lambda_u\quad\text{for every Borel set}\quad  B.
\end{eqnarray*}

$\textbf{IV}.$ Let us assume the contrary, i.e. $h$ has the representation \eqref{QID_fol_general_case}, however $\inf_{t\in \R^d} |h(t)| = 0$. Since $e^z \ne 0$ for all $z \in \CC$, then $h(t) \ne 0$ for all $t \in \R^d$. Hence it is sufficient to focus on the case when $h$ has the representation \eqref{QID_fol_general_case}, $h(t) \ne 0$ for all $t \in \R^d$, and $\inf_{t \in \R^d}|h(t)| = 0$.  

Due to \eqref{QID_fol_general_case}, for every fixed $\tau \in \R^d$ we have the following representation
\begin{eqnarray*}
	\frac{h(t+\tau)h(t-\tau)}{h^2(t)} = \exp\biggl\{-\tfrac{1}{2}\lng \tau,Q\tau\rng + 2\int\limits_{\R^d\setminus\{\oo\}} e^{i\lng t,u\rng} \bigr( \cos\bigr(\lng\tau,u\rng\bigr) - 1\bigr) \nu(\dd u)\biggr\},\quad t \in \R^d. 
\end{eqnarray*}
It follows that 
\begin{eqnarray*}
	\left| \frac{h(t+\tau)h(t-\tau)}{h^2(t)}\right| \leqslant \exp\Biggr\{\biggl( \tfrac{1}{2}\|Q\| + \int\limits_{0<\|u\|<1} \|u\|^2 \bigr|\nu\bigr|(\dd u) \biggr) \|\tau\|^2 + 4\int\limits_{\|u\|>1} \bigr|\nu\bigr|(\dd u)\Biggr\},\quad t \in \R^d.
\end{eqnarray*}
Hence, for every $\tau \in \R^d$ there exists $C_{\tau}$ such that 
\begin{eqnarray*}
	\sup_{t\in \R^d} \left| \frac{h(t+\tau)h(t-\tau)}{h^2(t)} \right|  \leqslant C_{\tau}. 
\end{eqnarray*}

Let $(t_n)_{n\in\N}$, $t_n \in \R^d$, be a sequence such that $h(t_n)$ tends to 0 as $n \to \infty$. If there exists $R > 0$ such that $\|t_n\| < R$ for every $n \in \N$, then there exists subsequence $(n_k)_{k \in \N}$  satisfying $t_{n_k} \to t_* \in \R^d$ as $k \to \infty$. Since $h$ is a continuous function,  $h(t_*) = 0$ that contradicts with the $(iv)$. It follows that $\|t_n\| \to \infty$ as $n \to \infty$.  Since $h$ is an almost periodic function, the sequence $( h(\cdot + t_n))_{n \in \N}$ is dense in the set of continuous functions, i.e. there exists a subsequence $(n_k)_{k \in \N}$ and a continuous function $\varphi$ such that 
\begin{eqnarray*}
	\sup_{\tau \in \R^d}\bigr| h(t_{n_k} + \tau) - \varphi(\tau)\bigr| \underset{k \to \infty} \longrightarrow 0.
\end{eqnarray*}
It is obvious that $|\varphi(\tau)| \leqslant C \defeq \sup_{t\in\R^d}|h(t)| < \infty$ for all $\tau \in \R^d$. Then 
\begin{eqnarray*}
	\Delta_k &\defeq& \sup_{\tau \in \R^d}\bigr| h(t_{n_k} + \tau)h(t_{n_k} - \tau) - \varphi(\tau)\varphi(-\tau)\bigr| \\ &\leqslant& \sup_{\tau \in \R^d}\bigr| h(t_{n_k} - \tau)\bigr| \cdot \bigr| h(t_{n_k} + \tau) - \varphi(\tau)\bigr| + \sup_{\tau \in \R^d} \bigr| h(t_{n_k} - \tau) - \varphi(-\tau)\bigr|  \cdot |\varphi(\tau)| \\ &\leqslant& 2C \sup_{\tau \in \R^d}\bigr| h(t_{n_k} + \tau) - \varphi(\tau)\bigr|  \underset{k \to \infty} \longrightarrow 0. 
\end{eqnarray*}
Let us assume that $\varphi(\tau)\varphi(-\tau) = 0$ for all $\tau \in \R^d$. It follows that 
\begin{eqnarray*}
	\sup_{\tau \in \R^d}\bigr| h(t_{n_k} + \tau)h(t_{n_k} - \tau)\bigr| \underset{k \to \infty} \longrightarrow 0. 
\end{eqnarray*}
So for any fixed $s \in \R^d$ 
\begin{eqnarray*}
	h(t_{n_k} + \tau)h(t_{n_k} - \tau)\Bigr|_{\tau = -t_{n_k} -s} = h(-s)h(2t_{n_k} +s) \underset{k \to \infty} \longrightarrow 0. 
\end{eqnarray*}
Since $h(s) \ne 0$ for every $s \in \R^d$, we have
\begin{eqnarray}\label{unif_limit}
	h(2t_{n_k} +s) \underset{k \to \infty} \longrightarrow 0. 
\end{eqnarray}
Next, it is easy to see that the function $h(2t_{n_k} + \cdot\,)$ is an almost periodic one. It means that there exists a subsequence $(n_{k_m})_{m \in \N}$ such that a sequence $\bigl( h(2t_{n_{k_m}} + \cdot\,)\bigr)_{m \in \N}$ has a uniform limit. From \eqref{unif_limit} one can conclude that
\begin{eqnarray*}
	\sup_{s\in \R^d} \bigr| h(2t_{n_{k_m}} + s)\bigr|\underset{m \to \infty} \longrightarrow 0. 
\end{eqnarray*}
Applying this with $s = -2t_{n_{k_m}}$, we come to a contradiction with $h(0)=1$. Therefore the assumption  $\inf_{t \in \R^d}|h(t)| = 0$ is false, i.e. $(i)$ follows from $(iv)$. 

$\textbf{V}.$ If $(i)$ holds, then $(ii)$ follows directly from \cite[Theorem 3.2]{BalanKrishtal}.\quad  $\Box$\\

\textbf{Proof of Theorem \ref{th_Repr}.}  The implication  $(a)\rightarrow(b)$ directly follows from the implication  $(iv)\rightarrow(i)$ of Theorem \ref{alper_mainres}.  The converse $(b)\rightarrow(a)$ holds due to  $(i)\rightarrow(iv)$ of Theorem \ref{alper_mainres} with applying \cite[Theorem 2.7]{BergKutLind} (so $\gamma \in \R^d$, $Q\in \R^{d\times d}$, $\nu$ is real-valued measure).  The representation \eqref{discr_Levy} holds due to $(iii)$ of Theorem \ref{alper_mainres} and \cite[Theorem 2.7]{BergKutLind} (so $\gamma \in \R^d$ and $\lambda_u\in \R$).\quad  $\Box$\\ 

\textbf{Proof of Theorem \ref{th:Cramer-Wold_QID}.} \textit{Necessity}. Due to  Theorem \ref{th_Repr}  and comments below, it is easily seen using formula \eqref{discr_Levy} that  if the distribution function  $F$ of a discrete random vector $\xi$ is (quasi-)infinitely divisible, then for any $c \in \R^d$ distribution functions $F_c$ of the random variables $\lng c, \xi\rng$ are (quasi-)infinitely divisible, respectively (there is the case $d=1$ for $F_c$). 

\textit{Sufficiency}. Let us consider a discrete random vector $\xi$ with distribution function  \eqref{main_rep} and characteristic function \eqref{ch_fun}. We write the latter in the expanded form: 
\begin{eqnarray*}
	f(t^{(1)},\ldots, t^{(d)}) = \sum_{k=1}^\infty p_{x_k}\exp\biggl\{i \sum_{j=1}^d t^{(j)}x^{(j)}_k \biggr\},
\end{eqnarray*}
where $x_k = ( x_k^{(1)}, \ldots, x_k^{(d)}) \in \R^d$ and $t^{(1)},\ldots, t^{(d)} \in \R$.

We now assume that the distribution functions $F_c$ of $\lng c, \xi\rng$  are quasi-infinitely divisible for any $c=(c^{(1)},\ldots, c^{(d)}) \in \R^d$. Let $f_c$ denote the corresponding characteristic functions. It is easily seen that
\begin{eqnarray*}
	f_c(t)=f(c^{(1)}t,\ldots, c^{(d)}t),\quad t\in \R.
\end{eqnarray*}
Applying Theorem \ref{th_Repr} to $F_c$ (here the case $d=1$), we conclude that there exists a constant $\mu_c>0$ such that
\begin{eqnarray}\label{assum_sepc0}
	\bigl|f(c^{(1)}t,\ldots, c^{(d)}t)\bigr|\geqslant \mu_c\quad \text{for all}\quad t \in \R.
\end{eqnarray}
In order to prove the quasi-infinite divisibility of $F$, according to  Theorem \ref{th_Repr}, it is sufficient to show that for some $\mu>0$
\begin{eqnarray} \label{cond_sepf0}
	\bigl|f(t^{(1)},\ldots, t^{(d)})\bigr| \geqslant \mu \quad \text{for all}\quad t^{(1)},\ldots, t^{(d)} \in \R.
\end{eqnarray}

We set  $X^{(j)} \defeq \{ x_k^{(j)}: p_{x_k}>0, k\in\N\} \subset \R$, $j = 1, \ldots, d$. Let us suppose that   $X^{(j)}\ne \{0\}$ for every $j=1, \ldots, d$, i.e. for every $j$ there exists $k\in\N$ such that $x_k^{(j)}\ne 0$. Therefore for every $j = 1, \ldots, d$ one can choose non-zero $\beta^{(j)}_l\in X^{(j)}$, $l\in \Index^{(j)}$ (here $\Index^{(j)}$ is at most countable index set) such that for every $k\in \N$ and for some  numbers $z^{(j)}_{k,l} \in \Z$ we have
\begin{eqnarray}\label{conc_xjk_betajl}
	x^{(j)}_k = \sum_{l\in \Index^{(j)}} z^{(j)}_{k,l} \beta^{(j)}_l,
\end{eqnarray}
where only finite number of $z^{(j)}_{k,l}$ are non-zero. Note that the numbers  $\beta^{(j)}_l$ can be chosen as linearly independent over $\Z$, that is the equation $z_1 \beta^{(j)}_{l_1} + \ldots+ z_{n}\beta^{(j)}_{l_n} = 0$ holds with $z_1, \ldots, z_{n} \in \Z$, and distinct $l_1, \ldots, l_n\in \Index^{(j)}$, $n\in\N$, if and only if $z_1 = \ldots= z_{n} = 0$. It follows that the numbers $z^{(j)}_{k,l}$ are uniquely determined for $x_k^{(j)}$ in \eqref{conc_xjk_betajl}. We observe that
\begin{eqnarray}\label{conc_Xj}
	\lng X^{(j)}\rng = \biggr\{ z_1 \beta^{(j)}_{l_1} + \ldots+ z_{n}\beta^{(j)}_{l_n}\colon \, z_1,\ldots, z_n \in \Z,\, l_1,\ldots, l_n\in\Index^{(j)},\, n\in\N\biggr\},\quad j=1,\ldots, d,
\end{eqnarray}
where $\lng\, \cdot\, \rng$ is determined by \eqref{modulus} (the one-dimesional case). 

We now propose the procedure of choosing of the numbers $c^{(1)},\ldots, c^{(d)}\in\R$  such that  the elements of the union system $\bigl\{c^{(1)} \beta^{(1)}_l: l\in \Index^{(1)}\bigr\}\cup\ldots \cup \bigl\{c^{(d)} \beta^{(d)}_l$: $l\in \Index^{(d)}\bigr\}$  are linearly independent over $\Z$. We first fix any $c^{(1)}\in \R\setminus \{0\}$. For every $v\in \lng X^{(2)}\rng\setminus \{0\}$ we define 
\begin{eqnarray*}
	D^{(2)}_v \defeq \bigl\{c\in\R:  c v \in \lng c^{(1)}X^{(1)}\rng \bigr\}.
\end{eqnarray*}
Here and below, for any set $X\subset \R$ we denote by $c X$, where $c\in\R$, the set $\{cx: x\in X\}$.   
Observe that every set $D^{(2)}_v$ is countable.  Then the set
\begin{eqnarray*}
	D^{(2)}\defeq \bigcup_{v\in \lng X^{(2)}\rng\setminus \{0\}} D^{(2)}_v
\end{eqnarray*}
is countable too. Hence the set  $C^{(2)}\defeq \R\setminus  D^{(2)}$  is not empty. We choose any $c^{(2)}\in C^{(2)}$. Observe that 
\begin{eqnarray*}
	C^{(2)}=  \R\setminus \bigcup_{v\in \lng X^{(2)}\rng\setminus \{0\}} D^{(2)}_v=\bigcap_{v\in \lng X^{(2)}\rng\setminus \{0\}}\R\setminus D^{(2)}_v.
\end{eqnarray*}
This means that for any $v\in \lng X^{(2)}\rng\setminus \{0\}$ the quantity $c^{(2)}v$  can not be a finite linear combination of elements $c^{(1)}\beta^{(1)}_{l}$, $l\in\Index^{(1)}$, with integer coefficients. Let $v=z_1 \beta^{(2)}_{l_1} + \ldots+ z_{n}\beta^{(2)}_{l_n}$ with some $z_1,\ldots, z_n \in \Z$.  $l_1,\ldots, l_n\in\Index^{(2)}$, and $n\in\N$. 
Since  $c^{(2)} v=z_1 (c^{(2)}\beta^{(2)}_{l_1}) + \ldots+ z_{n}(c^{(2)}\beta^{(2)}_{l_n})$, by the above argument, the elements in the union system $\bigl\{c^{(1)}\beta^{(1)}_l: l\in \Index^{(1)}\bigr\}\cup \bigl\{c^{(2)} \beta^{(2)}_l: l\in \Index^{(2)}\bigr\}$ are linear independent over $\Z$. We next consider the set of all finite linear combitations of $\bigl\{c^{(1)}\beta^{(1)}_l: l\in \Index^{(1)}\bigr\}\cup \bigl\{c^{(2)} \beta^{(2)}_l$: $l\in \Index^{(2)}\bigr\}$ with integer coefficients. It is the set $\lng c^{(1)}X^{(1)}\cup c^{(2)}X^{(2)} \rng$. For every $v\in \lng X^{(3)}\rng\setminus \{0\}$ we define 
\begin{eqnarray*}
	D^{(3)}_v \defeq \bigl\{c\in\R:  c v \in \lng c^{(1)}X^{(1)}\cup c^{(2)}X^{(2)} \rng \bigr\}.
\end{eqnarray*}
Every $D^{(3)}_v$ is countable. Hence the set
\begin{eqnarray*}
	D^{(3)}\defeq \bigcup_{v\in \lng X^{(3)}\rng\setminus \{0\}} D^{(3)}_v
\end{eqnarray*}
is  countable too. Since the set  $C^{(3)}\defeq \R\setminus  D^{(3)}$  is not empty, we choose any $c^{(3)}\in C^{(3)}$. Observe that 
\begin{eqnarray*}
	C^{(3)}=  \R\setminus \bigcup_{v\in \lng X^{(3)}\rng\setminus \{0\}} D^{(3)}_v=\bigcap_{v\in \lng X^{(3)}\rng\setminus \{0\}}\R\setminus D^{(3)}_v.
\end{eqnarray*}
Hence for any $v\in \lng X^{(3)}\rng\setminus \{0\}$ the quantity $c^{(3)}v$  can not be a finite linear combination of elements of $\bigl\{c^{(1)}\beta^{(1)}_l: l\in \Index^{(1)}\bigr\}\cup \bigl\{c^{(2)} \beta^{(2)}_l: l\in \Index^{(2)}\bigr\}$ with integer coefficients. This implies  that the elements in the union system $\bigl\{c^{(1)}\beta^{(1)}_l: l\in \Index^{(1)}\bigr\}\cup \bigl\{c^{(2)} \beta^{(2)}_l: l\in \Index^{(2)}\bigr\}\cup \bigl\{c^{(3)} \beta^{(3)}_l: l\in \Index^{(3)}\bigr\}$ are linear independent over $\Z$. We next proceed  analogously and thus we obtain that  the elements of the union system $\bigl\{c^{(1)} \beta^{(1)}_l: l\in \Index^{(1)}\bigr\}\cup\ldots \cup \bigl\{c^{(d)} \beta^{(d)}_l$: $l\in \Index^{(d)}\bigr\}$  are linearly independent over $\Z$ as required.

We now prove \eqref{cond_sepf0}. Suppose, contrary to our claim, that \eqref{cond_sepf0} is false, i.e.  for any $\e>0$ there exist   $t^{(1)}_\e,\ldots, t^{(d)}_\e\in\R$ such that $|f(t^{(1)}_\e,\ldots, t^{(d)}_\e)|\leqslant \e$. So we fix $\e>0$ and such $t^{(j)}_\e$, $j=1,\ldots, d$. We first find $N_\e\in\N$ such that $\sum_{k=N_\e+1}^\infty p_{x_k}\leqslant\e$ (see \eqref{main_rep} and \eqref{ch_fun},   $\sum_{k=1}^\infty p_{x_k}=1$,  $p_{x_k}\geqslant 0$). Then
\begin{eqnarray}\label{conc_sumN}
	\sup_{t^{(1)},\ldots, t^{(d)} \in \R}\biggl|\sum_{k=N_\e+1}^\infty p_{x_k} \exp\biggl\{i\sum_{j=1}^d x_k^{(j)}t^{(j)}\biggr\}\biggr|\leqslant\sum_{k=N_\e+1}^\infty p_{x_k}\leqslant\e.
\end{eqnarray}
Hence we get
\begin{eqnarray*}
	\bigl|f(t^{(1)}_\e,\ldots, t^{(d)}_\e)\bigr| &=& \biggl|\sum_{k\in\N} p_{x_k} \exp\biggl\{i\sum_{j=1}^d x_k^{(j)}t^{(j)}_\e\biggr\}\biggr|\\
	&\geqslant& \biggl|\sum_{k=1}^{N_\e} p_{x_k} \exp\biggl\{i\sum_{j=1}^d x_k^{(j)}t^{(j)}_\e\biggr\}\biggr|- \biggl|\sum_{k=N_\e+1}^\infty p_{x_k} \exp\biggl\{i\sum_{j=1}^d x_k^{(j)}t^{(j)}_\e\biggr\}\biggr|\\
	&\geqslant& \biggl|\sum_{k=1}^{N_\e} p_{x_k} \exp\biggl\{i\sum_{j=1}^d x_k^{(j)}t^{(j)}_\e\biggr\}\biggr|- \e.	
\end{eqnarray*}
Due to representations \eqref{conc_xjk_betajl}, we write: 
\begin{eqnarray}
	\sum_{k=1}^{N_\e} p_{x_k} \exp\biggl\{i\sum_{j=1}^d x_k^{(j)}t^{(j)}_\e\biggr\}&=&\sum_{k=1}^{N_\e} p_{x_k} \exp\biggl\{i\sum_{j=1}^d \biggl(\sum_{l\in \Index^{(j)}} z^{(j)}_{k,l} \beta^{(j)}_l\biggr)t^{(j)}_\e\biggr\}\nonumber\\
	&=&\sum_{k=1}^{N_\e} p_{x_k} \exp\biggl\{i\sum_{j=1}^d \sum_{l\in \Index^{(j)}} z^{(j)}_{k,l} \bigl(\beta^{(j)}_l t^{(j)}_\e\bigr)\biggr\}.\label{def_sumNe}
\end{eqnarray}
Let us fix $c^{(1)}, \ldots, c^{(d)} \in \R$ such that the elements of the union system $\bigl\{c^{(1)} \beta^{(1)}_l: l\in \Index^{(1)}\bigr\}\cup\ldots \cup \bigl\{c^{(d)} \beta^{(d)}_l: l\in \Index^{(d)}\bigr\}$  are linearly independent over $\Z$. By the Kronecker theorem (see \cite[p.37]{LevitanZhikov}), for any $\delta>0$ we can find $t'_\delta$ such that all following inequalities  hold with some integers $n_l^{(j)}$:
\begin{eqnarray}\label{Kroneckersystem}
	\bigl| c^{(j)}\beta^{(j)}_l t'_\delta -\beta^{(j)}_l t^{(j)}_\e-2\pi n_l^{(j)}\bigr|<\delta, \quad l\in\Index^{(j)}_\e,\quad j=1,\ldots, d,
\end{eqnarray}
where $\Index^{(j)}_\e$ is the set of all $l\in \Index^{(j)}$ such that $z^{(j)}_{k,l}\ne 0$ for some $k=1,\ldots, N_\e$. Since only finite number of $z^{(j)}_{k,l}$ are non-zero in \eqref{conc_xjk_betajl}, the set $\Index^{(j)}_\e$ is finite and the system \eqref{Kroneckersystem} has only finite number of inequalities. 
Let us choose $\delta=\delta_\e$ such that
\begin{eqnarray}\label{def_deltae}
	\delta_\e\cdot \max_{k=1,\ldots, N_\e}\biggl\{\sum_{j=1}^d \sum_{l\in \Index^{(j)}_\e} |z^{(j)}_{k,l}|\biggr\}\leqslant\e.
\end{eqnarray}
Observe that
\begin{eqnarray*}
	\Delta_\e&\defeq&\biggl|\sum_{k=1}^{N_\e} p_{x_k} \exp\biggl\{i\sum_{j=1}^d \sum_{l\in \Index^{(j)}_\e} z^{(j)}_{k,l} c^{(j)}\beta^{(j)}_l t'_{\delta_\e}\biggr\}-\sum_{k=1}^{N_\e} p_{x_k} \exp\biggl\{i\sum_{j=1}^d \sum_{l\in \Index^{(j)}_\e} z^{(j)}_{k,l} \beta^{(j)}_l t^{(j)}_\e\biggr\} \biggr|\\
	&\leqslant& \sum_{k=1}^{N_\e} p_{x_k} \biggl|\exp\biggl\{i\sum_{j=1}^d \sum_{l\in \Index^{(j)}_\e} z^{(j)}_{k,l} \bigl(c^{(j)}\beta^{(j)}_l t'_{\delta_\e}-\beta^{(j)}_l t^{(j)}_\e\bigr)\biggr\}-1\biggr|\\
	&=& \sum_{k=1}^{N_\e} p_{x_k} \biggl|\exp\biggl\{i\sum_{j=1}^d \sum_{l\in \Index^{(j)}_\e} z^{(j)}_{k,l} \bigl(c^{(j)}\beta^{(j)}_l t'_{\delta_\e}-\beta^{(j)}_l t^{(j)}_\e-2\pi n_l^{(j)}\bigr)\biggr\}-1\biggr|.
\end{eqnarray*}	
The last equality holds because all $z^{(j)}_{k,l}$ and  $n_l^{(j)}$  are integer. Next, using well known inequality $|e^{iy}-1|\leqslant |y|$, $y\in\R$, and applying \eqref{Kroneckersystem} and \eqref{def_deltae}, we obtain 	
\begin{eqnarray*}	
	\Delta_\e&\leqslant& \sum_{k=1}^{N_\e} \biggl(p_{x_k}\sum_{j=1}^d \sum_{l\in \Index^{(j)}_\e}\Bigl( \bigl|z^{(j)}_{k,l}\bigr|\cdot \bigl|c^{(j)}\beta^{(j)}_l t'_{\delta_\e}-\beta^{(j)}_l t^{(j)}_\e -2\pi n_l^{(j)}\bigr|\Bigr)\biggr)\\
	 &\leqslant&  \max_{k=1,\ldots, N_\e}\biggl\{\sum_{j=1}^d \sum_{l\in \Index^{(j)}_\e} |z^{(j)}_{k,l}|\cdot\delta_\e\biggr\}\cdot  \sum_{k=1}^{N_\e} p_{x_k} \\
	 &\leqslant&  \delta_\e\cdot\max_{k=1,\ldots, N_\e}\biggl\{\sum_{j=1}^d \sum_{l\in \Index^{(j)}_\e} |z^{(j)}_{k,l}|\biggr\}\leqslant\e.
\end{eqnarray*}
Returning to \eqref{def_sumNe}, we have
\begin{eqnarray*}
	\biggl|\sum_{k=1}^{N_\e} p_{x_k} \exp\biggl\{i\sum_{j=1}^d \sum_{l\in \Index^{(j)}} z^{(j)}_{k,l} \beta^{(j)}_l t^{(j)}_\e\biggr\} \biggr|\geqslant\biggl|\sum_{k=1}^{N_\e} p_{x_k} \exp\biggl\{i\sum_{j=1}^d \sum_{l\in \Index^{(j)}} z^{(j)}_{k,l} c^{(j)}\beta^{(j)}_l t'_{\delta_\e}\biggr\}\biggr|-\e.
\end{eqnarray*} 
Note that here  we write  $\Index^{(j)}$ instead of $\Index^{(j)}_\e$  that is obviously possible by the definition of $\Index^{(j)}_\e$.
Thus we get
 \begin{eqnarray*}
 	\bigl|f(t^{(1)}_\e,\ldots, t^{(d)}_\e)\bigr| 
 	&\geqslant& \biggl|\sum_{k=1}^{N_\e} p_{x_k} \exp\biggl\{i\sum_{j=1}^d \sum_{l\in \Index^{(j)}} z^{(j)}_{k,l} c^{(j)}\beta^{(j)}_l t'_{\delta_\e}\biggr\}\biggr|- 2\e.	
 \end{eqnarray*}
According to \eqref{conc_xjk_betajl}, we next write 
 \begin{eqnarray*}
 	\sum_{k=1}^{N_\e} p_{x_k} \exp\biggl\{i\sum_{j=1}^d \sum_{l\in \Index^{(j)}} z^{(j)}_{k,l} c^{(j)}\beta^{(j)}_l t'_{\delta_\e}\biggr\}=\sum_{k=1}^{N_\e} p_{x_k} \exp\biggl\{i\sum_{j=1}^d c^{(j)} x^{(j)}_{k}  t'_{\delta_\e}\biggr\}.
 \end{eqnarray*}
Due to \eqref{conc_sumN}, we get
\begin{eqnarray*}
	\biggl|\sum_{k=1}^{N_\e} p_{x_k} \exp\biggl\{i\sum_{j=1}^d c^{(j)} x^{(j)}_{k}  t'_{\delta_\e}\biggr\}\biggr|\geqslant \biggl|\sum_{k=1}^{\infty} p_{x_k} \exp\biggl\{i\sum_{j=1}^d c^{(j)} x^{(j)}_{k}  t'_{\delta_\e}\biggr\}\biggr|-\e=\bigl|f(c^{(1)}t'_{\delta_\e},\ldots, c^{(d)} t'_{\delta_\e})\bigr|-\e.
\end{eqnarray*}
So we have
\begin{eqnarray*}
	\e\geqslant|f(t^{(1)}_\e,\ldots, t^{(d)}_\e)| 
	\geqslant |f(c^{(1)}t'_{\delta_\e},\ldots, c^{(d)} t'_{\delta_\e})|- 3\e.	
\end{eqnarray*}
Thus for any $\e>0$ we found $t'_{\delta_\e}$ such that
\begin{eqnarray*}
	|f(c^{(1)}t'_{\delta_\e},\ldots, c^{(d)} t'_{\delta_\e})|\leqslant 4\e.
\end{eqnarray*}
This obviously contradicts to the assumption \eqref{assum_sepc0}.  So \eqref{cond_sepf0} holds. 

We proved the Cram\'er--Wold device for the quasi-infinite divisibility. Let us now consider the case of infinite divisibility. Let the distribution functions $F_c$ of random variables $\lng c, \xi \rng$ be infinitely divisible for any $c \in \R^d$. Then they are also quasi-infinitely divisible. From what has already been proved, the distribution function $F$ of the random vector $\xi$ is also quasi-infinitly divisible and, by Theorem \ref{th_Repr}, its characteristic function $f$ admits the representation
\begin{eqnarray*}
	f(t) = \exp\biggr\{ i\lng t,\gamma \rng + \sum_{u \in \langle X \rangle \setminus\{\oo\}} \lambda_u \bigr(e^{i\lng t,u\rng} - 1\bigr) \biggr\}, \quad t\in \R^d,
\end{eqnarray*}
where $\gamma \in \lng X \rng $, $\lambda_u \in \R$ for all $u \in \langle X \rangle \setminus\{\oo\}$, and $\sum_{u \in \langle X \rangle \setminus\{0\}}|\lambda_u| < \infty$. What is left is to show that $\lambda_u \geqslant 0$ for all $u \in \langle X \rangle \setminus\{\oo\}$. Let us write the characteristic function $f_c$ of $F_c$ for any $c \in \R^d$:
\begin{eqnarray*}
     f_c(t) = \exp\biggl\{it\lng c, \gamma\rng + \sum_{u \in \lng X \rng \setminus\{\oo\}}\lambda_u \bigr(e^{it\lng c, u \rng} - 1\bigr)\biggr\}, \quad t \in \R. 
\end{eqnarray*}
Let us fix $c = (c^{(1)}, \ldots, c^{(d)})\in \R^d$ such that the elements of the union system $\{c^{(1)} \beta^{(1)}_l: l\in \Index^{(1)}\}\cup\ldots \cup \{c^{(d)} \beta^{(d)}_l: l\in \Index^{(d)}\}$  are linearly independent over $\Z$. On account of \eqref{conc_xjk_betajl},  \eqref{conc_Xj}, and that $\lng X\rng= \lng X^{(1)} \rng\times \ldots \times \lng X^{(d)} \rng$,  we have $\lng c, u_1 \rng \ne \lng c, u_2 \rng$ for any distinct $u_1, u_2 \in \lng X \rng\setminus\{\oo\}$. Since $F_c$ is infinitely divisible (by assumption),  we conclude that $\lambda_u \geqslant 0$ for all $u \in \lng X \rng\setminus\{\oo\}$. \quad $\Box$

\section{Acknowledgments}
The work of A. A. Khartov was supported by RFBR--DFG grant 20-51-12004. The work of I. A. Alexeev was supported in part by the M\"obius Contest Foundation for Young Scientists.

\end{document}